\documentclass{elsart}

\usepackage{amsmath,amssymb}

\newcounter{sub}
\newenvironment{sub}%
{\begin{list}{(\arabic{sub})}{\usecounter{sub}%
\setlength{\leftmargin}{2em}}}{\end{list}}

\def\char{\mbox{char }}
\def\im{\mbox{Im }}

\begin{document}
\begin{frontmatter}

\title{Negative solutions to three-dimensional monomial Noether problem}
\author{Aiichi Yamasaki}
\address{Department of Mathematics, Kyoto University, Japan}

\begin{abstract}
Three-dimensional monomial Noether problem can have negative solutions
for 8 groups by the suitable choice of the coefficients.
We find the necessary and sufficient condition for the coefficients
to have a negative solution.
The results are obtained by two criteria of irrationality
using Galois cohomology.
\end{abstract}

\end{frontmatter}

\section{Introduction}

Let $G$ be a finite subgroup of $GL(n,\mathbb{Z})$
and $K$ be a field.

We define the action of $G$ on $K(x_1,\cdots,x_n)$
as follows.

For all $\sigma =(a_{ij})\ \in G$,
$$ x_i^\sigma= c_{\sigma,i} \prod _{j=1} ^n x_j^{a_{ij}} \quad (i=1,\ldots ,n),
\quad c_{\sigma,i} \in K^\times. $$

Let $K(x_1,\cdots,x_n)^G$ be the subfield of $G$-invariant functions of 
$K(x_1,\cdots,x_n)$.
{\it Is $K(x_1,\cdots,x_n)^G$ a rational extension of $K$?}
This is called the monomial Noether problem.
Especially when all $c_{\sigma,i}$'s are $1$,
the problem is called the purely monomial Noether problem.

In general, for mutually conjugate subgroups 
$G_1$ and $G_2$ of $GL(n,\mathbb{Z})$, monomial Noether problems  are the same by 
changing variables. 
So we have only to consider the conjugacy classes of 
finite subgroups of $GL(n,\mathbb{Z})$. 
Hence we have 73 kinds of problems in the case of $n=3$. See \cite{cryst4}.  

Hereafter we shall always assume $\char K \not=2$,
then 2 dimensional monomial Noether problems
and 3 dimensional purely monomial Noether problems
are all known to be affirmative \cite{Hajja1,Hajja2,Hajja3,Hajja4,Hoshi}.
In this paper, we shall consider 3 dimensional non-purely monomial Noether problems.
(Each of the above mentioned 73 kinds of problems consists of many problems
by the choice of coefficients $c_{\sigma,i}$). 

Among 73 conjugacy classes of finite subgroups of $GL(3,\mathbb{Z})$,
63 ones have affirmative solutions to monomial Noether problem.
This means that the answer is affirmative for any choice
of the coefficients $c_{\sigma,i}$.

8 ones have negative solutions,
which means that the answer is negative for some choice of the coefficients
$c_{\sigma,i}$, while it is affirmative for other choice of $c_{\sigma,i}$.

The remaining 2 cases are unsolved.
The answer is affirmative for some choice of $c_{\sigma,i}$,
including purely monomial ones,
but the answer is not known for other choice of $c_{\sigma,i}$.

In this paper, we shall discuss on the negative 8 cases.
The positive 63 cases and the unsolved 2 cases
will be discussed in another paper.

In order to prove the negativity of the Noether problem,
Galois cohomology is sometimes useful.
In \S\ref{cohomology} we shall give two criteria,
and call them the non-vanishing cohomology test and the parity test.
The idea of the former test was given by Lenstra \cite{Lenstra}
and used by many authors.
The idea of the latter test was given by Saltman \cite{Saltman}
to prove Theorem \ref{3121} mentioned later.

Among 8 negative cases, Theorem \ref{3422} and Theorem \ref{3121}
are already known cases,
obtained by Kang \cite[Theorem 1.8]{Kang} and Saltman \cite{Saltman} respectively.

Similar discussions can be applied to prove the negativity of other 6 cases.
So, for the sake of comparison,
we include the proof of Theorems \ref{3422} and \ref{3121} in full,
though they are not new results.

\section{Two criteria of Irrationality} \label{cohomology}

Galois cohomology is sometimes useful to prove the negativity of
the monomial Noether problem \cite{Lenstra,Saltman}.

Let $L$ be a finite Galois extension of $K$ with the Galois group $\mathfrak{G}$.
$\mathfrak{G}$ acts on $L(x_1,x_2,\cdots,x_n)$,
assuming that it acts on each $x_i$ trivially.
Similarly $G$ acts on $L(x_1,x_2,\cdots,x_n)$,
assuming that it acts on $L$ trivially.
Then we have

$$K(x_1,\cdots,x_n)^G=L(x_1,\cdots,x_n)^{\mathfrak{G}\times G}
=\big(L(x_1,\cdots,x_n)^G\big)^\mathfrak{G}.$$
Suppose that $L(x_1,\cdots,x_n)^G$ is rational over $L$
and $L(x_1,\cdots,x_n)^G=L(y_1,\cdots,y_n)$.
Then $K(x_1,\cdots,x_n)^G$ is rational over $K$,
if and only if $L(y_1,\cdots,y_n)^\mathfrak{G}$ is so.
Let $\mathcal{M}$ be a finite subset of irreducible polynomials in $L[y_1,\cdots,y_n]$,
which contains all $y_i$'s.
Let $M$ be the free Abelian group generated by $\mathcal{M}$.
It is a subgroup of $L(y_1,\cdots,y_n)^\times/L^\times$.
We assume that $M$ is closed under the action of $\mathfrak{G}$,
so that $M$ is a $\mathfrak{G}$-module.

Here we suppose that the transcendental basis $y_1,\cdots,y_n$ is explicitly known.

Assume that $L(y_1,\cdots,y_n)^\mathfrak{G}$ is rational and
$L(y_1,\cdots,y_n)^\mathfrak{G}=K(u_1,\cdots,u_n)$.
This time the assumption is virtual,
and assume only the existence of some $u_1,\cdots,u_n \in K(x_1,\cdots,x_n)$.
Since $L(y_1,\cdots,y_n)=L(u_1,\cdots,u_n)$,
the quotient fields of two rings $R_1=L[M]$ and $R_2=L[u_1,\cdots,u_n]$ coincide.
($R_1$ is explicitly determined while $R_2$ is virtual).
So that we have
$$\exists r_1 \in R_1^\mathfrak{G}, R_2 \subset R_1\big[\frac{1}{r_1}\big], \quad
\exists r_2 \in R_2^\mathfrak{G}, R_1 \subset R_2\big[\frac{1}{r_2}\big].$$
Put $R:=R_1\big[\frac{1}{r_1r_2}\big]=R_2\big[\frac{1}{r_1r_2}\big]$.
Then we have
$$R^\times=R_1^\times P=R_2^\times Q$$
where $P$ is the free module generated by irreducible factors
of $r_1r_2$ in $L[y_1,\cdots,y_n]$ other than those in $\mathcal{M}$,
and $Q$ is the free module generated by irreducible factors
of $r_1r_2$ in $L[u_1,\cdots,u_n]$.
($M$ is explicit, but $P$ and $Q$ are virtual).
Since $R_1^\times=L^\times M$, $R_2^\times=L^\times$,
we have $MP \simeq Q \pmod{L^\times}$.
All of $M$, $P$, $Q$ are $\mathfrak{G}$-modules.
Though $Q$ is the direct sum of $M$ and $P$ as $\mathbb{Z}$-modules,
it is not the direct sum as $\mathfrak{G}$-modules in general.

$P$ and $Q$ are permutation modules
(i.e. they are free $\mathbb{Z}$-modules
and $\mathfrak{G}$ acts as permutation of the basis),
and as such we have $H^1(\mathfrak{G},Q)=\hat H^{-1}(\mathfrak{G},Q)=0$,
where $\hat H$ is Tate cohomology.

$M$ is a direct sum factor of $Q$ if and only if $H^1(\mathfrak{H},M)=0$
for all subgroups $\mathfrak{H}$ of $\mathfrak{G}$.
The only if part is evident.
The reason of the if part is as follows.
We write the permutation of $\{p_j\}$ induced by $\sigma \in \mathfrak{G}$
by the same $\sigma$,
then we can write $p_j^\sigma=m_j(\sigma)+p_{\sigma(j)}$
for some $m_j(\sigma) \in M$.
(We write $MP$ additively though
it is a subgroup of $L(y_1,\cdots,y_n)^\times/L^\times$).
Fix any $j$ and let $\mathfrak{H}$ be the stabilizer of $j$,
namely $\mathfrak{H}=\{\sigma \in \mathfrak{G} \, | \, \sigma(j)=j\}$,
then $m_j(\sigma)$ is a $\mathfrak{H}$-cocycle,
so that if $H^1(\mathfrak{H},M)=0$
we have $\exists m_j \in M$,
$m_jp_j$ is $\mathfrak{H}$-invariant.

From this, we see that it is possible to modify the representative system
$\{p_j\}$ of $P=Q/M$ to make all $m_j(\sigma)$ zero,
which means that $M$ is a direct sum factor of $Q$.

If $M$ is a direct sum factor of $Q$,
we must also have $\hat H^{-1}(\mathfrak{G},M)=0$.
This proves the validity of the following criterion.

\underline{\bf The non-vanishing cohomology test}
\hfill\break
{\it If $\forall \mathfrak{H} < \mathfrak{G}, \, H^1(\mathfrak{H},M)=0$
and $\hat H^{-1}(\mathfrak{G},M) \not= 0$,
then $L(y_1,\cdots,y_n)^\mathfrak{G}$ is not rational over $K$.}

Even when $\hat H^{-1}(\mathfrak{G},M)=0$,
sometimes irrationality can be proved using another criterion.
The following is the parity test.

When $Q \simeq M \oplus P$,
we write the injection $M \rightarrow Q$ by $i$,
the projection $Q \rightarrow M$ by $p$,
then $p \circ i=id_M$.
Both $i$ and $p$ are $\mathfrak{G}$-homomorphisms,
and $\mathfrak{G}$ acts as a permutation module on $Q$
and acts in the explicitly given way on $M$.
Using this fact,
we examine the image of $y_j$ by $p \circ i$,
and check the parity (even or odd) of its coefficient.
If the coefficient of $y_j$ of this image is even,
then $p \circ i=id_M$ is not possible,
so that $L(y_1,\cdots,y_n)^\mathfrak{G}$ can not be rational over $K$.

Saltman proved Theorem \ref{3121} in \S\ref{s3121}
by the contradiction $p \circ i(y_j) \in 2M$.
But more generally, we get the following result.

\underline{\bf The parity test}
\hfill \break
{\it Assume that for any $\mathfrak{G}$-homomorphisms
$i: M \rightarrow Q$ and $p: Q \rightarrow M$,
we have $p \circ i(y_j) \not\equiv y_j \mod 2$
for some $j$,
then $L(y_1,\cdots,y_n)^\mathfrak{G}$ is not rational over $K$.}

\section{The problem $R(a,b,c)$.} \label{Rabc}

In this section, some preliminary discussions are given.

\begin{lem} \label{lem1}
The following two-dimensional Noether problems are affirmative.
Namely, both of $K(x,y)^{<\sigma>}$ and $K(x,y)^{<\tau>}$
are rational over $K$.

(1) $\sigma: x \mapsto -x, \,
y \mapsto \frac{ax^2+b}{y}, \, a,b \in K^\times$

(2) $\tau: x \mapsto \frac{a}{x}, \,
y \mapsto \frac{b(x+\frac{a}{x})+c}{y}, \,
a,b,c \in K^\times$.
\end{lem}

\begin{pf}
(1) We have $K(x,y)^{<\sigma>}=K(z_1,z_2,z_3)$
where $z_1=x^2$, $z_2=y+\frac{ax^2+b}{y}$,
$z_3=\big(y-\frac{ax^2+b}{y}\big)\big/x$.
Since $z_1z_3^2=z_2^2-4(az_1+b)$,
we have $z_1 \in K(z_2,z_3)$
so that $K(x,y)^{<\sigma>}=K(z_2,z_3)$.

(2) We have $K(x,y)^{<\tau>}=K(z_1,z_2,z_3)$
where $z_1=x+\frac{a}{x}$,
$z_2=y+\frac{b(x+\frac{a}{x})+c}{y}$,
$z_3=\big\{y-\frac{b(x+\frac{a}{x})+c}{y}\big\}\big/(x-\frac{a}{x})$.
Since $(z_1^2-4a)z_3^2=z_2^2-4(bz_1+c)$,
we get $z_2^2=d^2z_1^2+ez_1+f$ for some $d,e,f \in K(z_3)$,
so that $K(z_3)(z_1,z_2)$ is rational over $K(z_3)$.
\end{pf}

\underline{\bf Problem $R(a,b,c)$} in the title of this section
is the Noether problem
for the following $\sigma$.
$$\sigma: x_1 \mapsto -x_1, \,
x_2 \mapsto \frac{a}{x_2}, \,
x_3 \mapsto \frac{-bx_1^2+c}{x_3}, \,
a,b,c \in K^\times$$
Problem $R(a,b,c)$ is equivalent with the rationality problem
of the quadratic extension of $K(z_1,z_2,z_3)$ given by
$z_0^2=(z_1^2-a)(z_2^2-b)(z_3^2-c)$.

The reason is as follows.
Since $K(x_1,x_2,x_3)^{<\sigma>}=K(u_1,u_2,u_3,u_4)$ where
$u_1=(x_2-\frac{a}{x_2})/x_1$,
$u_2=x_2+\frac{a}{x_2}$,
$u_3=x_3+\frac{-bx_1^2+c}{x_3}$,
$u_4=\big(x_3-\frac{-bx_1^2+c}{x_3}\big)\big/(x_2-\frac{a}{x_2})$,
and since $(u_2^2-4a)u_4^2=u_3^2-4(-bx_1^2+c)$,
$x_1^2=\frac{u_2^2-4a}{u_1^2}$,
we have $u_1^2=(u_2^2-4a)(u_1^2u_4^2-4b)\big/(u_3^2-4c)$.
So if we put $z_1=u_2, z_2=u_1u_4, z_3=u_3$ and $z_0=(u_3^2-4c)u_1$,
the problem is reduced to the second form stated above.
Note that $R(a,b,c)$ depends only on the $K^{\times2}$-cosets
of $a,b,c$.

\begin{lem} \label{Rabclem}
Let $L=K(\sqrt{a},\sqrt{b},\sqrt{c})$.
$R(a,b,c)$ is affirmative
if and only if $[L:K]\leq 2$
or $[L:K]=4, abc \not\in K^{\times2}$.
In other words, it is affirmative if and only if
at least one of $a,b,c,ab,ac,bc$ belongs to $K^{\times2}$.
\end{lem}

\begin{pf}
$R(a,b,c)$ is affirmative
if $a=\alpha^2 \in K^{\times2}$,
because putting $z_0^\prime=\frac{z_0}{z_1+\alpha}$,
we have $z_0^{\prime2}=\frac{z_1-\alpha}{z_1+\alpha}(z_2^2-b)(z_3^2-c)$,
so that $z_1 \in K(z_0^\prime,z_2,z_3)$,
hence $K(z_0,z_1,z_2,z_3)=K(z_0^\prime,z_1,z_2,z_3)=K(z_0^\prime,z_2,z_3)$.
Similar discussion holds when $b$ or $c \in K^{\times2}$.

In the action of $\sigma$,
replace $x_3$ by $x_2x_3$,
then $b$ changes to $ab$,
and $c$ changes to $ac$.
Therefore $R(a,b,c)$ is equivalent to $R(a,ab,ac)$,
so that $R(a,b,c)$ is affirmative when $ab \in K^{\times2}$.
From the symmetricity  of the problem,
the same holds when $ac$ or $bc \in K^{\times2}$.

The case $[L:K]=8$ is easily reduced to
the case $[L:K]=4, \, abc \in K^{\times2}$ as follows.
Suppose $[L:K]=8$,
then we have $[L:K^\prime]=4$ for $K^\prime=K(\sqrt{abc})$
and $abc \in K^{\prime\times2}$,
so the negativity of $R(a,b,c)$ over $K^\prime$
implies the negativity over $K$.

We shall prove the negativity for the case
$[L:K]=4, \, abc \in K^{\times2}$.
We can assume $c=ab$.
Put $\alpha=\sqrt{a}$ and $\beta=\sqrt{b}$,
then $L=K(\alpha,\beta)$ and
$L(x_1,x_2,x_3)^{<\sigma>}$ is rational over $L$.

Putting $x_2^\prime=\frac{x_2-\alpha}{x_2+\alpha}$ and
$x_3^\prime=\frac{x_3-\beta x_1-\alpha\beta}{x_3+\beta x_1+\alpha\beta}$,
we have $\sigma: x_1 \mapsto -x_1$, $x_2^\prime \mapsto -x_2^\prime$,
$x_3^\prime \mapsto -x_3^\prime$ so that
$L(x_1,x_2,x_3)^{<\sigma>}=L(y_1,y_2,y_3)$ where
$y_1=x_1^2$, $y_2=x_1x_2^\prime$, $y_3=x_1x_3^\prime$.

$\mathfrak{G}=Gal(L/K)$ is isomorphic to $C_2 \times C_2$,
and generated by $\tau_1: \alpha \mapsto -\alpha, \, \beta \mapsto \beta$
and $\tau_2: \alpha \mapsto \alpha, \, \beta \mapsto -\beta$.
$\mathfrak{G}$ acts on $y_1,y_2,y_3$ etc. as follows.
{\small
\begin{center}
\begin{tabular}{c|ccccc}
 & $y_1$ & $y_2$ & $y_3$ & $y_1+\alpha y_3$ & $y_3+\alpha$ \\
\hline
$\tau_1$ & $y_1$ & $y_1/y_2$ & $(y_3+\alpha)y_1\big/(y_1+\alpha y_3)$ &
$y_1(y_1-a)\big/(y_1+\alpha y_3)$ &
$y_3(y_1-a)\big/(y_1+\alpha y_3)$ \\
$\tau_2$ & $y_1$ & $y_2$ & $y_1/y_3$ & $y_1(y_3+\alpha)\big/y_3$
& $(y_1+\alpha y_3)\big/y_3$ \\
\end{tabular}
\end{center}
}
Let $M$ be the $\mathbb{Z}$-module of rank 6 generated by
$y_1$,$y_2$,$y_3$,$y_1-a$,$y_1+\alpha y_3$ and $y_3+\alpha$.
Then the action of $\mathfrak{G}$ on $M$ is represented as matrices as follows.
$$m(\tau_1)=\begin{pmatrix}
 1 &   &   &   &   &   \\
 1 &-1 &   &   &   &   \\
 1 &   & 0 &   &-1 & 1 \\
   &   &   & 1 &   &   \\
 1 &   &   & 1 &-1 &   \\
   &   & 1 & 1 &-1 & 0
\end{pmatrix},
m(\tau_2)=\begin{pmatrix}
 1 &   &   &   &   &   \\
   & 1 &   &   &   &   \\
 1 &   &-1 &   &   &   \\
   &   &   & 1 &   &   \\
 1 &   &-1 &   & 0 & 1 \\
   &   &-1 &   & 1 & 0
\end{pmatrix}.$$

From this, we can calculate to get $H^1(\mathfrak{H},M)=0$
for any proper subgroup $\mathfrak{H}$ of $\mathfrak{G}$,
but $H^1(\mathfrak{G},M) \simeq \mathbb{Z}/2\mathbb{Z}$.
The non-trivial element of  $H^1(\mathfrak{G},M)$ is given by
$a(\tau_1)=1, \, a(\tau_2)=a(\tau_1\tau_2)=\frac{y_3(y_3+\alpha)}{y_1+\alpha y_3}$.

We construct a $\mathbb{Z}$-module $M^\prime$ of rank 7,
by adding $y_3(y_3+\alpha)+y_1+\alpha y_3=y_3^2+2\alpha y_3+y_1$
as the seventh generator.
Then the action of $\mathfrak{G}$ on $M^\prime$ is given by
$$m^\prime(\tau_1)=\begin{pmatrix}
m(\tau_1) & {\bf 0} \\
\begin{matrix} 1 & 0 & 0 & 1 &-2 & 0 \end{matrix} & 1
\end{pmatrix},
m^\prime(\tau_2)=\begin{pmatrix}
m(\tau_2) & {\bf 0} \\
\begin{matrix} 1 & 0 &-2 & 0 & 0 & 0 \end{matrix} & 1
\end{pmatrix},$$
and we have $H^1(\mathfrak{G},M^\prime)=0$.

As for $\hat{H}^{-1}$, a calculation shows that
$\hat{H}^{-1}(\mathfrak{G},M) \simeq \mathbb{Z}/2\mathbb{Z}$.
The non-trivial element is $y_2/y_3$,
which does not vanish in $\hat{H}^{-1}(\mathfrak{G},M^\prime)$,
thus $\hat{H}^{-1}(\mathfrak{G},M^\prime) \simeq \mathbb{Z}/2\mathbb{Z}$.
Therefore $K(x_1,x_2,x_3)^{<\sigma>}$ is not rational
by the non-vanishing cohomology test,
so $R(a,b,ab)$ is negative.
\end{pf}

\section{The group (3,4,2,2)}

The result here (Theorem \ref{3422}) is already known (Kang, \cite[Theorem 1.8]{Kang}),
but we shall give the proof in full as a standard example of the non-vanishing cohomology test.
(The numbering (3,4,2,2) of the group follows the GAP list \cite{cryst4}).
The group is isomorphic to $C_4$ and generated by the following $\sigma$.
$$C_4 \simeq <\sigma>, \quad
\sigma:x_1 \mapsto ax_2, \, x_2 \mapsto bx_3, \, x_3 \mapsto \frac{c}{x_1x_2x_3}.$$

Multiplying each $x_i$ by a constant factor, we can set $a=b=1$,
so the problem depends essentially only on $c$.
Then the problem is as follows:

``For what values of $c$ the field $K(x_1,x_2,x_3)^{<\sigma>}$
is rational over $K$?''
Here $\sigma$ is given by
$$\sigma:x_1 \mapsto x_2 \mapsto x_3 \mapsto \frac{c}{x_1x_2x_3} \mapsto x_1.$$

\begin{thm} \label{3422}
The monomial Noether problem for the group (3,4,2,2)
is affirmative if and only if $-1 \in K^{\times 2}$ or $c \in K^{\times2}$
or $c \in -4K^{\times4}$.
\end{thm}

\begin{pf}
$\sigma^2$ maps as $x_1 \leftrightarrow x_3$, $x_2 \mapsto \frac{c}{x_1x_2x_3}$
so that it keeps $x_1^\prime=x_1x_3$ invariant
and acts monomially on $x_2, x_3$
with $K(x_1^\prime)$-coefficient.
Thus we get $K(x_1,x_2,x_3)^{<\sigma^2>}=K(y_1,y_2,y_3)$
where $y_1=x_1^\prime=x_1x_3$,
$y_2=\frac{\xi_0}{\xi_2}$, $y_3=\frac{\xi_2}{\xi_3}$
with $\xi_2=x_2-\frac{c}{x_1x_2x_3}$,$\xi_3=x_3-x_1$,
$\xi_0=x_2x_3-\frac{c}{x_2x_3}$.

Now $\sigma$ maps as $y_1 \mapsto \frac{c}{y_1}$,
$y_3 \mapsto -\frac{1}{y_3}$, $y_2 \mapsto \frac{-y_1y_3+\frac{c}{y_1y_3}}{y_2}$
so that $y_3^\prime:=y_1y_3 \mapsto -\frac{c}{y_3^\prime}$
$y_2 \mapsto -\big(y_3^\prime-\frac{c}{y_3^\prime}\big)\big/y_2$.
If $c \in K^{\times2}$,
then $y_1^\prime :=\frac{y_1-\sqrt{c}}{y_1+\sqrt{c}}\big(y_3^\prime+\frac{c}{y_3^\prime}\big)$
is $\sigma$-invariant,
so $K(x_1,x_2,x_3)^{<\sigma>}=K(y_1,y_2,y_3)^{<\sigma>}=K(y_1^\prime)(y_2,y_3^\prime)^{<\sigma>}$ is rational
by Lemma \ref{lem1} (2).

If $-1 \in K^{\times2}$,
the same holds because
$y_1^{\prime\prime}:=\frac{y_3-\sqrt{-1}}{y_3+\sqrt{-1}}(y_3^\prime+\frac{c}{y_3^\prime})$
is $\sigma$-invariant.

If $-c=d^2 \in K^{\times2}$,
then $y_3^{\prime\prime}:=\frac{y_3^\prime-d}{y_3^\prime+d} \mapsto -y_3^{\prime\prime}$
and $y_2 \mapsto \frac{d}{y_2}\left(\frac{y_3^{\prime\prime}+1}{y_3^{\prime\prime}-1}
+\frac{y_3^{\prime\prime}-1}{y_3^{\prime\prime}+1}\right)
=\frac{2d(y_3^{\prime\prime2}+1)}{y_2(y_3^{\prime\prime2}-1)}$
so that $y_2^\prime:=(1+y_3^{\prime\prime})y_2
\mapsto \frac{-2d(y_3^{\prime\prime2}+1)}{y_2^\prime}$.
Together with $y_1 \mapsto \frac{-d^2}{y_1}$,
the Noether problem is reduced to $R(-1,2d,-2d)$,
so that it is affirmative if and only if one of $-1, \pm2d$
belongs to $K^{\times2}$ by Lemma \ref{Rabclem},
but $\pm2d \in K^{\times2}$ is equivalent with $c \in -4K^{\times4}$.

If none of $-1, \pm c$ belongs to $K^{\times2}$,
put $K^\prime=K(d)$ with $-c=d^2$,
then none of $-1, \pm2d$ belongs to $K^{\prime\times2}$
so that the Noether problem is negative over $K^\prime$,
hence negative over $K$.

This completes the proof of Theorem \ref{3422}.
\end{pf}

\section{The group (3,1,2,1)} \label{s3121}

A standard example of the application of the parity test is provided by the group (3,1,2,1).
The group is isomorphic to $C_2$, and generated by
$$\sigma: x_1 \mapsto \frac{a_1}{x_1}, \, x_2 \mapsto \frac{a_2}{x_2}, \,
x_3 \mapsto \frac{a_3}{x_3} \quad (a_1,a_2,a_3 \in K^\times).$$

The following result was obtained by Saltman \cite{Saltman},
but for the sake of comparison we include the proof in full.

Saltman seems to think that the key of the proof is the concept of retract rationality.
We admit the importance of retract rationality,
but think that the key of the proof is the parity test given in \S \ref{cohomology}.
Though we follow Saltman's proof essentially,
the stress is put on the application of the parity test
whose calculation is given in detail.

\begin{thm} \label{3121}
The monomial Noether problem for the group (3,1,2,1) is affirmative
if and only if $[K(\sqrt{a_1},\sqrt{a_2},\sqrt{a_3}):K] \leq 4$.
\end{thm}

\begin{pf}
Evidently the problem depends only on $K^{\times2}$-cosets of $a_1,a_2,a_3$.
Let $\sqrt{a_i}=\alpha_i$ and $L=K(\alpha_1,\alpha_2,\alpha_3)$.

If $a_1 \in K^{\times2}$,
then $x_1^\prime=\frac{x_1-\alpha}{x_1+\alpha}(x_2-\frac{a_2}{x_2})$ is $\sigma$-invariant,
and the problem is reduced to two-dimensional monomial Noether problem over $K(x_1^\prime)$,
so that it is affirmative.
The same holds when $a_2$ or $a_3 \in K^{\times2}$.
If we replace $x_1$ by $x_1x_2$,
we see that the problem is affirmative if $a_1a_2 \in K^{\times2}$.
The same holds when $a_1a_3$, $a_2a_3$ or $a_1a_2a_3 \in K^{\times2}$.
Thus, the problem is affirmative if $[L:K] \leq 4$.

Now we shall prove the negativity for the case $[L:K]=8$.
The proof is due to Saltman \cite{Saltman}.

$Gal(L/K)=\mathfrak{G}$ is isomorphic to $C_2 \times C_2 \times C_2$,
whose generators are $\tau_i \, (1 \leq i \leq 3)$
such that $\tau_i: \alpha_i \mapsto -\alpha_i, \,
\alpha_j \mapsto \alpha_j (j \not= i)$.
Add one more variable $x_0$ and assume that $\sigma$ acts as $x_0 \mapsto \frac{1}{x_0}$.
We shall prove that
$K(x_0,x_1,x_2,x_3)^{<\sigma>}$ is not rational.

$L(x_0,x_1,x_2,x_3)^{<\sigma>}$ is rational over $L$,
and the transcendental basis is given by
$y_0=x_0^{\prime2}$, $y_1=x_0^\prime x_1^\prime$,
$y_2=x_0^\prime x_2^\prime$, $y_3=x_0^\prime x_3^\prime$
where $x_0^\prime=\frac{x_0-1}{x_0+1}$,
$x_i^\prime=\frac{x_i-\alpha_i}{x_i+\alpha_i} \, (1 \leq i \leq 3)$.

Clearly $y_0$ is $\mathfrak{G}$-invariant,
and $\tau_i$ maps $x_i^\prime$ to $\frac{1}{x_i^\prime}$,
and $x_j^\prime$ to $x_j^\prime \, (j \not= i)$,
so that it maps $y_i$ to $\frac{y_0}{y_i}$,
and $y_j$ to $y_j \, (j \not= i)$.

Let $M$ be the $\mathbb{Z}$-module generated by $y_0,y_1,y_2$ and $y_3$.
The action of $\mathfrak{G}$ is represented by matrices as follows.
$$m(\tau_1)=\begin{pmatrix}
 1 & 0 & 0 & 0 \\
 1 &-1 & 0 & 0 \\
 0 & 0 & 1 & 0 \\
 0 & 0 & 0 & 1
\end{pmatrix}, \,
m(\tau_2)=\begin{pmatrix}
 1 & 0 & 0 & 0 \\
 0 & 1 & 0 & 0 \\
 1 & 0 &-1 & 0 \\
 0 & 0 & 0 & 1
\end{pmatrix}, \,
m(\tau_3)=\begin{pmatrix}
 1 & 0 & 0 & 0 \\
 0 & 1 & 0 & 0 \\
 0 & 0 & 1 & 0 \\
 1 & 0 & 0 &-1
\end{pmatrix}.$$
From this, we can calculate the cohomology groups as follows.

$\mathfrak{G}$ has 15 non-trivial subgroups including $\mathfrak{G}$ itself.

For each subgroup of order 2, we have
\hfill\break
$H^1(<\tau_i>,M)=0 \quad i=1,2,3$
\hfill\break
$H^1(<\tau_i\tau_j>,M) \simeq \mathbb{Z}/2\mathbb{Z} \quad (i \not= j)$
\hfill\break
$H^1(<\tau_1\tau_2\tau_3>,M) \simeq \left(\mathbb{Z}/2\mathbb{Z}\right)^2$
\hfill\break
The non-trivial element of $H^1(<\tau_1\tau_2>,M)$ is given by
the 1-cocycle $a(\tau_1\tau_2)=y_0^{-1}y_1y_2$.
\hfill\break
The non-trivial elements of $H^1(<\tau_1\tau_2\tau_3>,M)$ are
the following $a_i \, (1 \leq i \leq 3)$;
$a_1(\tau_1\tau_2\tau_3)=y_0^{-1}y_1y_2$,
$a_2(\tau_1\tau_2\tau_3)=y_0^{-1}y_1y_3$,
$a_3(\tau_1\tau_2\tau_3)=y_0^{-1}y_2y_3$.

For subgroups of order 4, we have
\hfill\break
$H^1(<\tau_i,\tau_j>,M)=0 \quad (i \not= j)$
\hfill\break
$H^1(<\tau_i,\tau_j\tau_k>,M) \simeq \mathbb{Z}/2\mathbb{Z}$ \quad
($i,j,k$ are mutually different)
\hfill\break
$H^1(<\tau_1\tau_2,\tau_1\tau_3>,M) \simeq \mathbb{Z}/2\mathbb{Z}$
\hfill\break
The non-trivial element of $H^1(<\tau_1,\tau_2\tau_3>,M)$
is given by $a(\tau_1)=1$, $a(\tau_2\tau_3)=a(\tau_1\tau_2\tau_3)=y_0^{-1}y_2y_3$.
\hfill\break
The non-trivial element of $H^1(<\tau_1\tau_2,\tau_1\tau_3>,M)$
is given by $a(\tau_1\tau_2)=y_0^{-1}y_1y_2$,
$a(\tau_1\tau_3)=y_0^{-1}y_1y_3$, $a(\tau_2\tau_3)=y_0^{-1}y_2y_3$.
\hfill\break
Finally, we have $H^1(\mathfrak{G},M)=0$.
\hfill\break
As for $\hat{H}^{-1}$, we have
$\hat{H}^{-1}(\mathfrak{G},M) \simeq \left(\mathbb{Z}/2\mathbb{Z}\right)^2$,
where non-trivial elements are $y_0^{-1}y_1y_2$, $y_0^{-1}y_1y_3$ and $y_0^{-1}y_2y_3$.

In order to make $H^1$ zero, we extend the $\mathbb{Z}$-module $M$ to a larger $M^\prime$.
Let $y_4=y_1+y_2$, $y_5=y_1y_2+y_0$, $y_6=y_1+y_3$, $y_7=y_1y_3+y_0$,
$y_8=y_2+y_3$, $y_9=y_2y_3+y_0$ and let $M^\prime$ be the $\mathbb{Z}$-module
of rank 10 generated by $y_0 \sim y_9$.

The action of $\mathfrak{G}$ on $M^\prime$ is represented by matrices as follows.
$$m^\prime(\tau_1)=\begin{pmatrix}
m(\tau_1) & & 0 & \\
A_1 & C & & \\
A_1 & & C & \\
0 & & & 1
\end{pmatrix}, \,
m^\prime(\tau_2)=\begin{pmatrix}
m(\tau_2) & & 0 & \\
A_2 & C & & \\
0 & & 1 & \\
A_2 & & & C
\end{pmatrix}, \,
m^\prime(\tau_3)=\begin{pmatrix}
m(\tau_3) & & 0 & \\
0 & 1 & & \\
A_3 & & C & \\
A_3 & & & C
\end{pmatrix}$$
where
$A_1=\begin{pmatrix}
 0 &-1 & 0 & 0 \\
 1 &-1 & 0 & 0
\end{pmatrix}$,
$A_2=\begin{pmatrix}
 0 & 0 &-1 & 0 \\
 1 & 0 &-1 & 0
\end{pmatrix}$,
$A_3=\begin{pmatrix}
 0 & 0 & 0 &-1 \\
 1 & 0 & 0 &-1
\end{pmatrix}$,
$C=\begin{pmatrix}
 0 & 1 \\
 1 & 0
\end{pmatrix}$.

The extended module $M^\prime$ makes $H^1$ zero for almost all $\mathfrak{H} < \mathfrak{G}$.
More strictly speaking,
except only one group $\mathfrak{H}_0=<\tau_1\tau_2,\tau_1\tau_3>$,
we have $H^1(\mathfrak{H},M^\prime)=0$ for $\mathfrak{H} \not= \mathfrak{H}_0$.

To make $H^1(\mathfrak{H}_0,M^{\prime\prime})=0$,
we add $y_{10}=y_1y_2+y_1y_3+y_2y_3+y_0$ and $y_{11}=y_0(y_1+y_2+y_3)+y_1y_2y_3$.
Then the action of $\mathfrak{G}$ on $M^{\prime\prime}$ is represented by
$$m^{\prime\prime}(\tau_i)=\begin{pmatrix}
m^\prime(\tau_i) & 0 \\
\begin{matrix} A_i & 0 \end{matrix} & C
\end{pmatrix}$$
This time we have $H^1(\mathfrak{H},M^{\prime\prime})=0$ for all $\mathfrak{H}<\mathfrak{G}$.
However the obstruction for $\hat{H}^{-1}(\mathfrak{G},M)$ is of the same kind
of that for $H^1(\mathfrak{H},M)$,
so $\hat{H}^{-1}$ also becomes zero by the extension of $M$.
We have $\hat{H}^{-1}(\mathfrak{G},M^\prime)=\hat{H}^{-1}(\mathfrak{G},M^{\prime\prime})=0$.
So the non-vanishing cohomology test can not be applied.

We shall follow the discussions in \S\ref{cohomology},
(using $M^{\prime\prime}$ instead of $M$).
$M^{\prime\prime}$ is a direct product factor of a permutation module $Q$,
so that we have $id_{M^{\prime\prime}}=p \circ i$,
where $i$ is the injection $M^{\prime\prime} \rightarrow Q$
and $p$ is the projection $Q \rightarrow M^{\prime\prime}$.
If we show that the image of $y_0$ by $p \circ i$ belongs to $2M^{\prime\prime}$,
a contradiction occurs
and $L(y_0,y_1,y_2,y_3)^\mathfrak{G}$ can not be rational over $K$.

First,we shall examine the image $i(y_0)$.
Since $y_0$ is $\mathfrak{G}$-invariant,
$i(y_0)$ is also $\mathfrak{G}$-invariant
and the exponent is constant on every transitive part of permutations by $\mathfrak{G}$.
Consider any transitive part $X$,
then $i(y_0)$ takes the form $m\sum_{q \in X} q$ on $X$.
(We write $Q$ additively).
Suppose that the stabilizer of $X$ is trivial,
hence $X$ consists of 8 irreducible factors.
Let $\beta=mq_0$ for a fixed $q_0 \in X$,
then $i(y_0)=\sum_{\tau \in \mathfrak{G}} \beta^\tau$ on $X$.

Suppose that the stabilizer of $X$ is not trivial.
Let $\tau$ be a non-trivial element of the stabilizer,
then at least one of $\alpha_i$ moves by $\tau$.
For simplicity, assume that $\alpha_1$ moves by $\tau$.
Since $\tau$ maps $y_1$ to $y_0y_1^{-1}$,
and since $q \in X$ does not move by $\tau$,
we have $n=m-n$ where $m$ and $n$ are the exponents of the factor $q$
of $i(y_0)$ and $i(y_1)$ respectively.
This implies $m=2n$, so $m$ is even.
The same discussion holds for $\alpha_2$ and $\alpha_3$ instead of $\alpha_1$,
so whenever the stabilizer of $X$ is not trivial,
the exponent $m$ should be even.
If $i(y_0)=m\sum_{q \in X}q$ on $X$,
then put $\gamma=\frac{m}{2}\left(\prod_{q \in X}q\right)$
and we have $i(y_0)=2\gamma$ on $X$.
Combining these two results, we get
$i(y_0)=\sum_{\tau \in \mathfrak{G}} \beta^\tau + 2\gamma$
for some $\beta,\gamma \in Q$.

The image by $p$ is written as
$$p \circ i(y_0)=\sum_{\tau \in \mathfrak{G}} p(\beta)^\tau + 2p(\gamma).$$
Evidently $2p(\gamma) \in 2M^{\prime\prime}$.
On the other hand,
the action of $\mathfrak{G}$ on $M^{\prime\prime}$ is represented by matrices
$m^{\prime\prime}(\tau)$,
so that $\sum_{\tau \in \mathfrak{G}} p(\beta)^\tau$
belongs to the image of $\sum_{\tau \in \mathfrak{G}}m^{\prime\prime}(\tau)$.

So, if every matrix elements of
$\sum_{\tau \in \mathfrak{G}} m^{\prime\prime}(\tau)$ is even,
then $\sum_{\tau \in \mathfrak{G}} p(\beta)^\tau$ belongs to $2M^{\prime\prime}$,
so that $p \circ i(y_0) \in 2M^{\prime\prime}$,
which is a contradiction.

The calculation of $\sum_{\tau \in \mathfrak{G}} m^{\prime\prime}(\tau)$
is easily done as follows.
$$\sum_{\tau \in \mathfrak{G}} m^{\prime\prime}(\tau)=\begin{pmatrix}
{\displaystyle \sum_{\tau \in \mathfrak{G}} m^\prime(\tau)} & {\bf 0} \\
\begin{matrix} E & \,\, & 0 \end{matrix} & F
\end{pmatrix}, \quad
E=\begin{pmatrix}
 4 &-4 &-4 &-4 \\
 8 &-4 &-4 &-4
\end{pmatrix}, \,
F=\begin{pmatrix}
 4 & 4 \\
 4 & 4
\end{pmatrix}$$
$$\sum_{\tau \in \mathfrak{G}} m^\prime(\tau)=\begin{pmatrix}
{\displaystyle \sum_{\tau \in \mathfrak{G}} m(\tau)} & & {\bf 0} & \\
G_1 & F & & \\
G_2 & & F & \\
G_3 & & & F
\end{pmatrix}, \quad
G_1=\begin{pmatrix}
 2 &-4 &-4 & 0 \\
 6 &-4 &-4 & 0
\end{pmatrix},$$
$$G_2=\begin{pmatrix}
 2 &-4 & 0 &-4 \\
 6 &-4 & 0 &-4
\end{pmatrix}, \,
G_3=\begin{pmatrix}
 2 & 0 &-4 &-4 \\
 6 & 0 &-4 &-4
\end{pmatrix}$$
$$\sum_{\tau \in \mathfrak{G}}m(\tau)=\begin{pmatrix}
8 & & & \\
4 & 0 & & \\
4 & & 0 & \\
4 & & & 0
\end{pmatrix}$$
Thus every entries of the matrix $\sum_{\tau \in \mathfrak{G}}m^{\prime\prime}(\tau)$
is even.
Therefore the negativity of the problem for $[L:K]=8$ has been proved,
which completes the proof of Theorem \ref{3121}.
\end{pf}

\section{The group (3,4,2,1)}

The group (3,4,2,1) is isomorphic to $C_4$, and generated by
$$\sigma: x_1 \mapsto \frac{a}{x_2}, \, x_2 \mapsto bx_1, \, x_3 \mapsto \frac{c}{x_3}.$$
By a suitable change of variables,
we can set $b=1$.
First we shall determine the fixed field of $\sigma^2$.
$$\sigma^2: x_1 \mapsto \frac{a}{x_1}, \, x_2 \mapsto \frac{a}{x_2}, \,
x_3 \mapsto x_3.$$
It is known that $K(x_1,x_2,x_3)^{<\sigma^2>}=K(y_1,y_2,y_3)$
where $y_1=(x_1+x_2)/(x_1x_2+a)$, $y_2=(x_1-x_2)/(x_1x_2-a)$, $y_3=x_3$.
The group $G=<\sigma>$ is of order 2 on $K(y_1,y_2,y_3)$
and acts as follows.
$$\sigma: y_1 \mapsto \frac{1}{ay_1}, \, y_2 \mapsto -\frac{1}{ay_2}, \,
y_3 \mapsto \frac{c}{y_3}.$$
This action is the same as that in (3,1,2,1),
so that we get
\begin{thm}
For the group (3,4,2,1) (with $b=1$),
the Noether problem is affirmative if and only if
$[K(\sqrt{a},\sqrt{-1},\sqrt{c}):K] \leq 4$.
\end{thm}

\section{The group (3,2,3,1)}

The group (3,2,3,1) is isomorphic to $C_2 \times C_2$,
and generated by the following $\sigma_1$ and $\sigma_2$.
$$\begin{cases}
\sigma_1: x_1 \mapsto \varepsilon_1 x_1, \,
x_2 \mapsto \varepsilon_2 x_2, \, x_3 \mapsto \frac{c}{x_3}
& (\varepsilon_1,\varepsilon_2=\pm1, \, c \in K^\times) \\
\sigma_2: x_1 \mapsto \frac{a}{x_1}, \, x_2 \mapsto \frac{b}{x_2}, \,
x_3 \mapsto \varepsilon_3 x_3 &
(\varepsilon_3=\pm1, \, a,b \in K^\times)
\end{cases}$$

\begin{thm} \label{3231thm}
The Noether problem for the group (3,2,3,1) is affirmative
if $\varepsilon_1=\varepsilon_2=1$ or $\varepsilon_3=1$,
and is reduced to $R(a,b,c)$ in \S\ref{Rabc}
if $\varepsilon_1=\varepsilon_2=\varepsilon_3=-1$.
\end{thm}

\begin{pf}
When $\varepsilon_3=1$, we have
$K(x_1,x_2,x_3)^{<\sigma_2>}=K(y_1,y_2,y_3)$ where
$$y_1=\big(x_1-\frac{a}{x_1}\big)\big/\big(x_1x_2-\frac{ab}{x_1x_2}\big), \,
y_2=\big(x_2-\frac{b}{x_2}\big)\big/\big(x_1x_2-\frac{ab}{x_1x_2}\big), \,
y_3=x_3,$$
and $\sigma_1$ acts on $y_i$ as
$$\sigma_1: y_1 \mapsto \varepsilon_2 y_1, \,
y_2 \mapsto \varepsilon_1 y_2, \, y_3 \mapsto \frac{c}{y_3},$$
so that $K(y_1,y_2,y_3)^{<\sigma_1>}=K(x_1,x_2,x_3)^{<\sigma_1,\sigma_2>}$ is rational.

When $\varepsilon_1=\varepsilon_2=1$,
we have $K(x_1,x_2,x_3)^{<\sigma_1>}=K(y_1,y_2,y_3)$ where
$y_1=x_1$, $y_2=x_2$, $y_3=x_3+\frac{c}{y_3}$,
and $\sigma_2$ acts on $y_i$ as
$y_1 \mapsto \frac{a}{y_1}$, $y_2 \mapsto \frac{b}{y_2}$,
$y_3 \mapsto \varepsilon_3 y_3$,
so that $K(y_1,y_2,y_3)^{<\sigma_2>}=K(x_1,x_2,x_3)^{<\sigma_1,\sigma_2>}$ is rational.

If either of $\varepsilon_1$ or $\varepsilon_2$ is $-1$,
replacing $x_1$ or $x_2$ by $x_1x_2$,
we can set $\varepsilon_1=\varepsilon_2=-1$.
Suppose that $\varepsilon_3=-1$.
Then we have $K(x_1,x_2,x_3)^{<\sigma_1>}=K(y_1,y_2,y_3)$
where $y_1=x_1x_2$, $y_2=\big(x_3-\frac{c}{x_3}\big)x_2$, $y_3=x_3+\frac{c}{x_3}$,
and $\sigma_2$ acts as
$y_1 \mapsto \frac{ab}{y_1}$, $y_2 \mapsto -\frac{b(y_3^2-4c)}{y_2}$, $y_3 \mapsto -y_3$.
Therefore the Noether problem is reduced to $R(ab,b,bc)$
which is equivalent with $R(a,b,c)$.

This completes the proof of Theorem \ref{3231thm}
\end{pf}

\section{The group (3,3,1,1)}

The group (3,3,1,1) is also isomorphic to $C_2 \times C_2$,
and generated by the following $\sigma_1$ and $\sigma_2$.
$$\begin{cases}
\sigma_1: x_1 \mapsto \varepsilon_1 x_1, \, x_2 \mapsto \frac{b}{x_2}, \,
x_3 \mapsto \frac{c}{x_3} &
(\varepsilon_1=\pm 1, \, b,c \in K^\times) \\
\sigma_2: x_1 \mapsto \frac{a}{x_1}, \, x_2 \mapsto \varepsilon_2 x_2, \,
x_3 \mapsto \frac{\varepsilon_3 c}{x_3} &
(\varepsilon_2, \varepsilon_3 = \pm 1, \, a \in K^\times)
\end{cases}$$
Note that $\sigma_1 \sigma_2: x_1 \mapsto \frac{\varepsilon_1 a}{x_1}, \,
x_2 \mapsto \frac{\varepsilon_2 b}{x_2}, \, x_3 \mapsto \varepsilon_3 x_3$.
The Noether problem depends only on the $K^{\times2}$-cosets of $a,b,c$.

\begin{thm} \label{3311thm}
\begin{sub}
\item[(1)] When $\varepsilon_1=\varepsilon_2=\varepsilon_3=1$,
the Noether problem is reduced to $R(a,b,c)$.
\item[(2)] When $\varepsilon_1=\varepsilon_2=-1, \, \varepsilon_3=1$,
the Noether problem is reduced to $R(a,b,c)$.
When $\varepsilon_1=\varepsilon_3=-1, \, \varepsilon_2=1$,
it is reduced to $R(-a,b,c)$,
and when $\varepsilon_2=\varepsilon_3=-1, \, \varepsilon_1=1$,
it is reduced to $R(a,-b,-c)$.
\item[(3)]
When $\varepsilon_1=\varepsilon_2=1, \, \varepsilon_3=-1$ or
$\varepsilon_1=\varepsilon_3=1, \, \varepsilon_2=-1$ or
$\varepsilon_2=\varepsilon_3=1, \, \varepsilon_1=-1$,
the Noether problem is always affirmative.
\item[(4)]
When $\varepsilon_1=\varepsilon_2=\varepsilon_3=-1$,
the results are as follows.
\begin{sub}
\item[{\romannumeral 1})]
When $a=\pm1$ or $b=\pm1$ or $c=\pm1$,
the Noether problem is reduced to $R(b,-1,c)$ or $R(a,-1,-c)$ or $R(a,-1,b)$
respectively.
\item[{\romannumeral 2})]
When $ab=1$ or $ac=-1$ or $bc=1$,
the Noether problem is reduced to $R(a,c,-ac)$ or $R(-a,b,ab)$ or $R(-b,a,ab)$
respectively.
\item[{\romannumeral 3})]
When none of $\pm a$, $\pm b$, $\pm c$, $ab$, $-ac$, $bc$ belongs to $K^{\times2}$,
the Noether problem is negative.
\end{sub}
\end{sub}
\end{thm}

\begin{pf}
\begin{sub}
\item
Suppose that $\varepsilon_1=\varepsilon_2=\varepsilon_3=1$,
then we have $K(x_1,x_2,x_3)^{<\sigma_1,\sigma_2>}=K(z_1,z_2,z_3,z_4)$
where $z_1=x_1+\frac{a}{z_1}$, $z_2=x_2+\frac{b}{x_2}$, $z_3=x_3+\frac{c}{x_3}$
and $z_4=\big(x_1-\frac{a}{x_1}\big)\big(x_2-\frac{b}{x_2}\big)
\big(x_3-\frac{c}{x_3}\big)$.

Since $z_4^2=(z_1^2-4a)(z_2^2-4b)(z_3^2-4c)$,
$K(x_1,x_2,x_3)^{<\sigma_1,\sigma_2>}$ is a quadratic extension of $K(z_1,z_2,z_3)$
with the above defining relation,
so that the Noether problem is reduced to $R(a,b,c)$.
\item
We shall consider the case $\varepsilon_1=\varepsilon_2=-1, \varepsilon_3=1$,
from which other two cases are derived by the symmetry of the problem.

Suppose that $\varepsilon_1=\varepsilon_2=-1, \varepsilon_3=1$,
then we have $K(x_1,x_2,x_3)^{<\sigma_1,\sigma_2>}=K(z_1,z_2,z_3,z_4)$
where $z_1=x_3+\frac{c}{x_3}$, $z_2=\big(x_1-\frac{a}{x_1}\big)\big(x_3-\frac{c}{x_3}\big)$,
$z_3=\big(x_1-\frac{a}{x_1}\big)\big(x_1+\frac{a}{x_1}\big)^{-1}\big(x_2+\frac{b}{x_2}\big)$
and $z_4=\big(x_2-\frac{b}{x_2}\big)\big(x_2+\frac{b}{x_2}\big)^{-1}\big(x_1+\frac{a}{x_1}\big)$.

Since we have $z_2^2=(z_1^2-4c)(z_4^2-4a)z_3^2(z_3^2-4b)^{-1}$,
the Noether problem is reduced to $R(a,b,c)$ by the same reason as (1).
\item
We shall consider the case $\varepsilon_1=\varepsilon_2=1, \varepsilon_3=-1$,
then we have $K(x_1,x_2,x_3)^{<\sigma_1,\sigma_2>}=K(z_1,z_2,z_3,z_4)$
where $z_1=x_1+\frac{a}{x_1}$,$z_2=x_2+\frac{b}{x_2}$,
$z_3=\big(x_1-\frac{a}{x_1}\big)^{-1}\big(x_3+\frac{c}{x_3}\big)$ and
$z_4=\big(x_2-\frac{b}{x_2}\big)^{-1}\big(x_3-\frac{c}{x_3}\big)$.

Since we have $z_1^2=\frac{z_4^2z_2^2}{z_3^2}+\frac{4}{z_3^2}(c-bz_4^2)+4a$,
regarding this identity as a relation of $z_1$ and $z_2$ over $K(z_3,z_4)$,
we see that the quadratic extension defined by this relation is rational.
\item
Suppose that $\varepsilon_1=\varepsilon_2=\varepsilon_3=-1$.
First assume $c=1$ (in general $c \in K^{\times2}$).

Putting $x_3^\prime=\frac{x_3-1}{x_3+1}$,
we have $\sigma_1: x_3^\prime \mapsto -x_3^\prime$,
$\sigma_2: x_3^\prime \mapsto -\frac{1}{x_3^\prime}$,
so that
$$\begin{cases}
\sigma_1:  x_1 \mapsto -x_1, \, x_3^\prime \mapsto -x_3^\prime, \,
x_2 \mapsto \frac{b}{x_2} & \\
\sigma_2: x_1 \mapsto \frac{a}{x_1}, \,
x_3^\prime \mapsto -\frac{1}{x_3^\prime}, \, x_2 \mapsto -x_2
\end{cases}$$
which is the same as the action of the group (3,2,3,1).
So by Theorem \ref{3231thm},
the Noether problem is reduced to $R(a,-1,b)$.

Starting from $c=-1$, we get the same result.
If $a=\pm1$ or $b=\pm1$, then the desired result in i)
is obtained by the symmetry of the problem.

Assume that $b=c$ (which is equivalent to $bc=1 \mod K^{\times2}$).
Then $\sigma_1$ and $\sigma_2$ act on $x_3^\prime=x_2x_3$ as
$\sigma_1: x_3^\prime \mapsto \frac{b^2}{x_3^\prime}$,
$\sigma_2: x_3^\prime \mapsto \frac{bx_2^2}{x_3^\prime}$,
so that we have $K(x_1,x_2,x_3)^{<\sigma_1>}=K(y_1,y_2,y_3)$
where $y_1=x_1\big(x_2-\frac{b}{x_2}\big)$, $y_2=x_2+\frac{b}{x_2}$
and $y_3=x_3^{\prime\prime}\big(x_2-\frac{b}{x_2}\big)$
with $x_3^{\prime\prime}=\frac{x_3^\prime-b}{x_3^\prime+b}$.

The action of $\sigma_2$ on $y_i$ is as follows.
$$\begin{cases}
y_1 \mapsto -\frac{a}{x_1}\big(x_2-\frac{b}{x_2}\big)
=-\frac{a}{y_1^2}(y_2^2-4b), & \\
y_2 \mapsto -y_2, & \\
y_3^\prime:=y_3-y_2=-\frac{2b(x_2+\frac{x_3^\prime}{x_2})}{x_3^\prime+b}
\mapsto \frac{2(x_2+\frac{bx_2}{x_3^\prime})}{\frac{x_2^2}{x_3^\prime}+1}
=\frac{-4b}{y_3^\prime} &
\end{cases}$$
Therefore $K(x_1,x_2,x_3)^{<\sigma_1,\sigma_2>}=K(y_1,y_2,y_3)^{<\sigma_2>}$
implies that the Noether problem is reduced to $R(-b,a,ab)$.

If $ab=1$ or $ac=-1$,
the desired result in ii) is obtained by the symmetry of the problem.

Suppose that none of $\pm a, \pm b, \pm c, ab, -ac, bc$
belongs to $K^{\times2}$.
Let $L=K(\sqrt{a},\sqrt{b},\sqrt{c},\sqrt{-1})$.
If $[L:K]=16$, then putting $K^\prime=K(\sqrt{c})$,
we have $c \in K^{\prime\times2}$
and $[L:K^\prime]=8$,
so the problem is negative over $K^\prime$,
hence negative over $K$.

If $[L:K]=8$, then one of $-ab, ac, -bc, \pm abc, -1$
belongs to $K^{\times2}$.
If $-ab$ or $-abc \in K^{\times2}$,
then $[L:K^\prime]=4$ and $-ab \in K^{\prime\times2}$,
so the problem is negative over $K^\prime$,
hence negative over $K$.
Similar discussions for $\pm a, \pm b, -c$ instead of $c$
assures the negativity of the problem
when $ac, -bc$ or $abc \in K^{\times2}$.

Putting $K^{\prime\prime}=K(\sqrt{bc})$,
we have $bc \in K^{\prime\prime\times2}$ and $[L:K^{\prime\prime}]=4$,
so that if $-1 \in K^{\times2}$,
then the problem is negative over $K^{\prime\prime}$,
hence over $K$.

If $[L:K]=4$, then three of $-ab, ac, -bc, \pm abc, -1$
belong to $K^{\times2}$.
The following two cases are possible.
\begin{sub}
\item[(A)]
$-1$ and $\pm abc$ belong to $K^{\times2}$.
\item[(B)]
$-ab, ac, -bc$ belong to $K^{\times2}$.
\end{sub}
We shall prove the negativity of the problem for (A) and (B)
by the non-vanishing cohomology test.

Put $\alpha=\sqrt{a}, \beta=\sqrt{b}, \gamma=\sqrt{c}$
and $x_2^\prime=\frac{x_2-\beta}{x_2+\beta}$,
$x_3^\prime=\frac{x_3-\gamma}{x_3+\gamma}$.
Then the action of $\sigma_1$ is
$\sigma_1: x_1 \mapsto -x_1, x_2^\prime \mapsto -x_2^\prime,
x_3^\prime \mapsto -x_3^\prime$,
therefore $L(x_1,x_2,x_3)^{<\sigma_1>}=L(y_1,y_2,y_3)$
where $y_1=x_1^2$, $y_2=x_1x_2^\prime$, $y_3=x_1x_3^\prime$.

The action of $\sigma_2$ is
$y_1 \mapsto \frac{a^2}{y_1}, y_2 \mapsto \frac{a}{y_2},
y_3 \mapsto -\frac{a}{y_3}$,
so that we have $L(x_1,x_2,x_3)^{<\sigma_1,\sigma_2>}=L(z_1,z_2,z_3)$
where $z_1=y_1^{\prime2}$, $z_2=y_1^\prime y_2^\prime$, $z_3=y_1^\prime y_3^\prime$
with
$y_1^\prime=\frac{y_1-a}{y_1+a}$, $y_2^\prime=\frac{y_2-\alpha}{y_2+\alpha}$,
$y_3^\prime=\frac{y_3-i\alpha}{y_3+i\alpha}$.

First, we consider the case (B).

$\mathfrak{G}=Gal(L/K)$ is isomorphic to $C_2 \times C_2$,
and generated by
$\tau_1: \alpha \mapsto -\alpha, i \mapsto i$ and
$\tau_2: \alpha \mapsto \alpha, i \mapsto -i$.
Evidently $z_1$ is $\mathfrak{G}$-invariant.
$\tau_1$ maps $y_2$ to $\frac{y_1}{y_2}$
and $y_3$ to $\frac{y_1}{y_3}$,
so that it maps $z_2$ to $z_1(z_2-1)(z_2-z_1)^{-1}$
and $z_3$ to $(z_3-z_1)(z_3-1)^{-1}$.
$\tau_2$ maps $y_2$ to $\frac{y_1}{y_2}$ and $y_3$ to $y_3$,
so that $z_2$ to $(z_2-z_1)(z_2-1)^{-1}$
and $z_3$ to $\frac{z_1}{z_3}$.

Let $M$ be the $\mathbb{Z}$-module of rank 8 generated by
$z_1$, $z_1-1$, $z_2$, $z_2-1$, $z_2-z_1$, $z_3$, $z_3-1$, $z_3-z_1$.
Then the action of $\mathfrak{G}$ is represented by the following matrices.
$$m(\tau_1)=\begin{pmatrix}
I_2 & 0 & 0 \\
B_1 & C_1 & 0 \\
B_2 & 0 & C_2
\end{pmatrix}, 
m(\tau_2)=\begin{pmatrix}
I_2 & 0 & 0 \\
B_2 & C_2 & 0 \\
B_3 & 0 & C_3
\end{pmatrix}, 
m(\tau_1\tau_2)=\begin{pmatrix}
I_2 & 0 & 0 \\
B_3 & C_3 & 0 \\
B_1 & 0 & C_1
\end{pmatrix}, $$
$$\begin{matrix}
B_1=\begin{pmatrix}
1 & \\ & 1 \\ 1 & 1
\end{pmatrix}, &
C_1=\begin{pmatrix}
 & 1 & -1 \\ 1 & & -1 \\ & & -1
\end{pmatrix}, \\
B_2=\begin{pmatrix}
0 & \\ & 1 \\ & 1
\end{pmatrix}, &
C_2=\begin{pmatrix}
 & -1 & 1 \\ & -1 & \\ 1 & -1 &
\end{pmatrix}, \\
B_3=\begin{pmatrix}
1 & \\ & 0 \\ 1 &
\end{pmatrix}, &
C_3=\begin{pmatrix}
-1 & & \\ -1 & & 1 \\ -1 & 1 &
\end{pmatrix}.
\end{matrix}$$

Next, we consider the case (A).
$\mathfrak{G}=Gal(L/K)$ is also isomorphic to $C_2 \times C_2$
and generated by $\tau_1: \beta \mapsto -\beta, \alpha \mapsto \alpha$
and $\tau_2:\beta \mapsto \beta, \alpha \mapsto -\alpha$.
Then we have
$$m(\tau_1)=\begin{pmatrix}
I_2 & 0 & 0 \\
B_2 & C_2 & 0 \\
B_1 & 0 & C_1
\end{pmatrix}, 
m(\tau_2)=\begin{pmatrix}
I_2 & 0 & 0 \\
B_3 & C_3 & 0 \\
B_2 & 0 & C_2
\end{pmatrix}, 
m(\tau_1\tau_2)=\begin{pmatrix}
I_2 & 0 & 0 \\
B_1 & C_1 & 0 \\
B_3 & 0 & C_3
\end{pmatrix}. $$
Interchanging $z_2$ and $z_3$,
the matrices become same with those for (B).
Namely, the action of $\mathfrak{G}$ is essentially same for (A) and (B),
so the Noether problem is same for both cases.

We shall return to the case (B).
We can calculate $H^1$ and $\hat{H}^{-1}$ as follows.
$H^1(<\tau_1>,M)=H^1(<\tau_2>,M)=H^1(<\tau_1\tau_2>,M)=0$,
$H^1(\mathfrak{G},M) \simeq \big(\mathbb{Z}/2\mathbb{Z}\big)^2$,
which is generated by
$a_1(\tau_1)=a_1(\tau_2)=\frac{z_2(z_1-1)}{(z_2-1)(z_2-z_1)}, a_1(\tau_1\tau_2)=1$
and
$a_2(\tau_1)=a_2(\tau_2)=\frac{z_3(z_3-1)}{z_3-z_1}, a_2(\tau_1\tau_2)=1$.
$\hat{H}^{-1}(\mathfrak{G},M) \simeq \big(\mathbb{Z}/2\mathbb{Z}\big)^2$,
which is generated by $\frac{z_3}{z_2}$ and $\frac{z_3-1}{z_2-1}$.

In order to make $H^1$ zero, we shall construct a $\mathbb{Z}$-module
$M^\prime$ of rank 10 by adding the following two generators.
$$z_2(z_1-1)+(z_2-1)(z_2-z_1)=(z_2-1)^2+z_1-1, \quad
z_3(z_3-1)+z_3-z_1=z_3^2-z_1.$$
Then the action of $\mathfrak{G}$ on $M^\prime$ is represented as
$$\begin{matrix}
m^\prime(\tau)=\begin{pmatrix}
m(\tau) & {\bf 0} \\ F(\tau) & I_2
\end{pmatrix}, &
F(\tau_1)=\begin{pmatrix}
 1 & 1 & 0 & 0 &-2 & 0 & 0 & 0 \\
 0 & 1 & 0 & 0 & 0 & 0 &-2 & 0
\end{pmatrix}, \\
& F(\tau_2)=\begin{pmatrix}
 0 & 1 & 0 &-2 & 0 & 0 & 0 & 0 \\
 1 & 0 & 0 & 0 & 0 &-2 & 0 & 0
\end{pmatrix}, \\
& F(\tau_1\tau_2)=\begin{pmatrix}
 1 & 0 &-2 & 0 & 0 & 0 & 0 & 0 \\
 1 & 1 & 0 & 0 & 0 & 0 & 0 &-2
\end{pmatrix}.
\end{matrix}$$
For this $M^\prime$, we get $H^1(\mathfrak{G},M^\prime)=0$,
but $\hat{H}^{-1}$ does not vanish at all,
so $\hat{H}^{-1}(\mathfrak{G},M^\prime)=\hat{H}^{-1}(\mathfrak{G},M)
\simeq \big(\mathbb{Z}/2\mathbb{Z}\big)^2$.
Therefore, by the non-vanishing cohomology test,
we see that the Noether problem is negative.
\end{sub}

The proof of Theorem \ref{3311thm} is now complete.
\end{pf}

\section{The group (3,4,3,1)}

The group (3,4,3,1) is isomorphic to $C_4 \times C_2$,
and generated by
$$\begin{cases}
\sigma_1: x_1 \mapsto \frac{a}{x_2}, \, x_2 \mapsto bx_1, \,
x_3 \mapsto \frac{c}{x_3} & (a,b,c \in K^\times) \\
\sigma_2: x_1 \mapsto \varepsilon_1 x_1, \, x_2 \mapsto \varepsilon_1 x_2, \,
x_3 \mapsto \frac{\varepsilon_2c}{x_3} & (\varepsilon_1,\varepsilon_2=\pm1)
\end{cases}$$
(If we consider $\sigma_2$ alone,
the general form is $x_1 \mapsto \varepsilon_1 x_1$,
$x_2 \mapsto \varepsilon_2 x_2$, $x_3 \mapsto \frac{d}{x_3}$,
but the condition that $\sigma_1$ and $\sigma_2$ should commute
implies $\varepsilon_1=\varepsilon_2$ and $\big(\frac{d}{c}\big)^2=1$).

By a suitable change of variables, we can set $b=1$.
Then $\sigma_1^2$ acts as
$x_1 \mapsto \frac{a}{x_1}$, $x_2 \mapsto \frac{a}{x_2}$, $x_3 \mapsto x_3$,
so that we have $K(x_1,x_2,x_3)^{<\sigma_1^2>}=K(y_1,y_2,y_3)$
where $y_1=\frac{x_1+x_2}{x_1x_2+a}$, $y_2=\frac{x_1-x_2}{x_1x_2-a}$,
$y_3=x_3$.
The actions of $\sigma_2$ and $\sigma_1\sigma_2$ on $y_i$ are as follows.
$$\begin{cases}
\sigma_2: y_1 \mapsto \varepsilon_1 y_1, \, y_2 \mapsto \varepsilon_1 y_2, \,
y_3 \mapsto \frac{\varepsilon_2 c}{y_3} & \\
\sigma_1 \sigma_2: y_1 \mapsto \frac{\varepsilon_1}{ay_1}, \,
y_2 \mapsto -\frac{\varepsilon_1}{ay_2}, \, y_3 \mapsto \varepsilon_2 y_3 &
\end{cases}$$
This action is the same as that in (3,2,3,1).
So we get the following result.

\begin{thm}
If $\varepsilon_1$ or $\varepsilon_2=1$,
the Noether problem of the group (3,4,3,1) (with $b=1$) is affirmative.
If $\varepsilon_1=\varepsilon_2=-1$,
then it is reduced to $R(a,-a,-c)$,
hence it is affirmative if and only if at least one of
$-1$, $\pm a$, $\pm ac$, $-c$ belongs to $K^{\times2}$.
\end{thm}

\section{The group (3,4,4,1)}

The group (3,4,4,1) is isomorphic to $D_4$,
and generated by
$$\begin{cases}
\sigma_1: x_1 \mapsto \frac{a}{x_2}, \, x_2 \mapsto bx_1, \,
x_3 \mapsto \alpha x_3 & (\alpha^4=1, \, a,b \in K^\times) \\
\sigma_2: x_1 \mapsto \varepsilon x_1, \, x_2 \mapsto \frac{\varepsilon ab}{x_2}, \,
x_3 \mapsto \frac{c}{x_3} & (\varepsilon=\pm1, \, c \in K^\times)
\end{cases}$$
(If we consider $\sigma_2$ alone,
the general form of $\sigma_2$ should be $x_2 \mapsto \frac{d}{x_2}$,
but the condition $\sigma_1\sigma_2=\sigma_2\sigma_1^{-1}$
implies $\frac{a}{d}=\frac{\varepsilon}{b}$,
hence $d=\varepsilon ab$).
By a suitable change of variables, we can set $b=1$.

\begin{thm} \label{3441thm}
If $-1 \in K^{\times2}$ or $\alpha=1$, then the Noether problem
of the group (3,4,4,1) (with $b=1$) is affirmative.
If $-1 \not\in K^{\times2}$ and $\alpha=-1$, then it is reduced to
$R(-1,-\varepsilon a,-c)$,
so that it is affirmative if and only if at least one of
$\pm a$, $\pm c$, $\varepsilon ac$ belongs to $K^{\times2}$.
\end{thm}

\begin{pf}
If $\alpha=\pm 1$, then  we have
$K(x_1,x_2,x_3)^{<\sigma_1^2>}=K(y_1,y_2,y_3)$ where
$y_1=\frac{x_1+x_2}{x_1x_2+a}$, $y_2=\frac{x_1-x_2}{x_1x_2-a}$, $y_3=x_3$.
The actions of $\sigma_1$ and $\sigma_1\sigma_2$ on $y_i$ are as follows.
$$\begin{cases}
\sigma_1: y_1 \mapsto \frac{1}{ay_1}, \, y_2 \mapsto -\frac{1}{ay_2}, \,
y_3 \mapsto \alpha y_3 & \\
\sigma_1\sigma_2: y_1 \mapsto \varepsilon y_1, \, y_2 \mapsto -\varepsilon y_2, \,
y_3 \mapsto \frac{\alpha c}{y_3} &
\end{cases}$$

This action is the same as that of (3,2,3,1).
So by Theorem \ref{3231thm},
the Noether problem is affirmative if $\alpha=1$,
and is reduced to $R(-1,-\varepsilon a, -c)$ if $\alpha=-1$.

If $\alpha^2=-1$, then $-1 \in K^{\times2}$.
We shall prove the affirmativity of the problem for this case.
Then $K(x_1,x_2,x_3)^{<\sigma_1^2>}=K(y_1,y_2,y_3^\prime)$
where $y_3^\prime=(x_1-\frac{a}{x_1})x_3$.
The action of $\sigma_1$ and $\sigma_1\sigma_2$ is
$$\begin{cases}
\sigma_1: y_3^\prime \mapsto
-\big(x_2-\frac{a}{x_2}\big)\alpha x_3=-\alpha\frac{y_1-y_2}{y_1+y_2}y_3^\prime & \\
\sigma_1\sigma_2: y_3^\prime \mapsto
\varepsilon \big(x_2-\frac{a}{x_2}\big)\frac{\alpha c}{x_3}
=\frac{4\varepsilon \alpha c}{y_3^\prime}
\frac{(1-ay_1^2)(1-ay_2^2)}{y_1^2-y_2^2} & 
\end{cases}$$
Then $\sigma_1\sigma_2$ maps $\eta_3=(y_1-y_2)y_3^\prime$ to
$\frac{4\alpha c}{\eta_3}(1-ay_1^2)(1-ay_2^2)$.
Suppose $\varepsilon=1$,
then $K(x_1,x_2,x_3)^{<\sigma_1^2,\sigma_1\sigma_2>}=K(z_1,z_2,z_3)$
where $z_1=y_1$, $z_2=\big(\eta_3-\eta_3^{\sigma_1\sigma_2}\big)\big/y_2$,
$z_3=\eta_3+\eta_3^{\sigma_1\sigma_2}$.

The action of $\sigma_1$ is $z_1 \mapsto \frac{1}{az_1}$,
$\eta_3 \mapsto \frac{-\alpha}{a}\big(\frac{1}{y_1}+\frac{1}{y_2}\big)
\frac{y_1-y_2}{y_1+y_2}y_3^\prime
=\frac{-\alpha}{ay_1y_2}\eta_3$,
so that $z_3 \mapsto \frac{-\alpha}{a}\frac{z_2}{z_1}$,
$z_2 \mapsto \frac{\alpha z_3}{z_1}$.
The action of $\sigma_1$ on $z_1, z_2, z_3$ is monomial,
and is not conjugate to (3,1,2,1),
which is the only negative group isomorphic to $C_2$.
So $K(z_1,z_2,z_3)^{<\sigma_1>}$ is rational.
(When $\varepsilon=-1$, by interchanging $y_1$ with $y_2$ and $a$ with $-a$,
the problem is reduced to the case of $\varepsilon=1$).
Thus, Theorem \ref{3441thm} has been proved for $-1 \in K^{\times2}$.
\end{pf}

\section{The group (3,3,3,1)}

The group (3,3,3,1) is isomorphic to $C_2 \times C_2 \times C_2$,
and generated by the following $\sigma_1$, $\sigma_2$ and $\sigma_3$.

$$G:\begin{cases}
\sigma_1: x_1 \mapsto \frac{a}{x_1}, \, x_2 \mapsto \varepsilon_{12} x_2, \,
x_3 \mapsto \varepsilon_{13} x_3 & \\
\sigma_2: x_1 \mapsto \varepsilon_{21} x_1, \, x_2 \mapsto \frac{b}{x_2}, \,
x_3 \mapsto \varepsilon_{23} x_3 & \\
\sigma_3: x_1 \mapsto \varepsilon_{31} x_1, \, x_2 \mapsto \varepsilon_{32} x_2, \,
x_3 \mapsto \frac{c}{x_3} & (a,b,c \in K^\times, \, \varepsilon_{ij} = \pm 1)
\end{cases}$$

\begin{thm} \label{3331}
The Noether problem of the group (3,3,3,1) is affirmative
except the following four cases.
\begin{sub}
\item[{\romannumeral 1})]
$\varepsilon_{12}=\varepsilon_{21}=1$ or
$\varepsilon_{13}=\varepsilon_{31}=1$ or
$\varepsilon_{23}=\varepsilon_{32}=1$,
all other $\varepsilon_{ij}$ are $-1$.
\item[{\romannumeral 2})]
all $\varepsilon_{ij}=-1$.
\end{sub}
For the exceptional four cases, it is affirmative if and only if
$[K(\sqrt{a},\sqrt{b},\sqrt{c}):K] \leq 4$.
\end{thm}

\begin{pf}
\begin{sub}
\item
When all $\varepsilon_{ij}$ are 1, $K(x_1,x_2,x_3)^G$ is rational.
The transcendental basis is given by
$z_1=x_1+\frac{a}{x_1}$, $z_2=x_2+\frac{b}{x_2}$, $z_3=x_3+\frac{c}{x_3}$.
\item
When $\varepsilon_{12}=\varepsilon_{13}=1$,
$K(x_1,x_2,x_3)^G$ is rational,
because $K(x_1,x_2,x_3)^{<\sigma_1>}=K(y_1,y_2,y_3)$
where $y_1=x_1+\frac{a}{x_1}$, $y_2=x_2$, $y_3=x_3$,
on which $\sigma_2$ and $\sigma_3$ act as
$$\begin{cases}
\sigma_2: y_1 \mapsto \varepsilon_{21} y_1, \, y_2 \mapsto \frac{b}{y_2}, \,
y_3 \mapsto \varepsilon_{23} y_3 & \\
\sigma_3: y_1 \mapsto \varepsilon_{31} y_1, \, y_2 \mapsto \varepsilon_{32} y_2, \,
y_3 \mapsto \frac{c}{y_3} &
\end{cases}$$
This action is reduced to a two-dimensional monomial one,
so that the Noether problem is affirmative.

The same holds when $\varepsilon_{21}=\varepsilon_{23}=1$
or $\varepsilon_{31}=\varepsilon_{32}=1$.
Among $2^6=64$ cases for the choice of $\varepsilon_{ij}$,
27 cases remain to be considered.
\item \label{2D2}
When $\varepsilon_{21}=\varepsilon_{31}=1$,
$K(x_1,x_2,x_3)^G$ is rational as shown below.

$K(x_1,x_2,x_3)^{<\sigma_1>}=K(y_1,y_2,y_3)$
where $y_1=x_1+\frac{a}{x_1}$ and for $i=2,3$,
$y_i=x_i$ if $\varepsilon_{1i}=1$,
and $y_i=\big(x_1-\frac{a}{x_1}\big)x_i$ if $\varepsilon_{1i}=-1$.
If $\varepsilon_{12}=-1$, then $\sigma_2$ maps $y_2$ to
$\varepsilon_{21}\frac{b(y_1^2-4a)}{y_2}$
and $\sigma_3$ maps $y_2$ to $\varepsilon_{31}\varepsilon_{32}y_2$.

Suppose that $\varepsilon_{21}=\varepsilon_{31}=1$,
then $y_1$ is $G$-invariant and the actions of $\sigma_2$ and $\sigma_3$
on $y_2$,$y_3$ are two-dimensional monomial over $K(y_1)$,
so that the Noether problem is affirmative.
The same holds when $\varepsilon_{12}=\varepsilon_{32}=1$
or $\varepsilon_{13}=\varepsilon_{23}=1$.
Thus 9 cases are settled, and 18 cases remain.
\item \label{3minus}
Suppose that $\varepsilon_{12}=\varepsilon_{23}=\varepsilon_{31}=1$,
$\varepsilon_{13}=\varepsilon_{21}=\varepsilon_{32}=-1$.
Then $\sigma_2$ and $\sigma_3$ act as
$$\begin{cases}
\sigma_2: y_1 \mapsto -y_1, \, y_2 \mapsto \frac{b}{y_2}, \,
y_3 \mapsto -y_3 & \\
\sigma_3: y_1 \mapsto y_1, \, y_2 \mapsto -y_2, \,
y_3 \mapsto \frac{c(y_1^2-4a)}{y_3}. &
\end{cases}$$
So $K(x_1,x_2,x_3)^{<\sigma_1,\sigma_2>}
=K(y_1,y_2,y_3)^{<\sigma_2>}=K(z_1,z_2,z_3)$,
where $z_1=\big(y_2-\frac{b}{y_2}\big)y_1$,
$z_2=y_2+\frac{b}{y_2}$, $z_3=\big(y_2-\frac{b}{y_2}\big)y_3$,
on which $\sigma_3$ acts as
$z_1 \mapsto -z_1, \, z_2 \mapsto -z_2, \,
z_3 \mapsto -\big(y_2-\frac{b}{y_2}\big)\frac{c(y_1^2-4a)}{y_3}
=-\frac{c}{z_3}\big\{z_1^2-4a(z_2^2-4b)\big\}$.

However, the Noether problem of $\sigma: x \mapsto -x, \, y \mapsto -y, \,
z \mapsto \frac{ax^2+by^2+c}{z}$ is always affirmative,
because $\xi=\frac{x}{y}$ is $\sigma$-invariant and
$z \mapsto \frac{(a\xi^2+b)y^2+c}{z}$,
whose numerator is a linear polynomial of $y^2$ over $K(\xi)$.

Thus in this case, $K(x_1,x_2,x_3)^G$ is rational.
The same holds also for $\varepsilon_{13}=\varepsilon_{21}=\varepsilon_{32}=1$,
$\varepsilon_{12}=\varepsilon_{23}=\varepsilon_{31}=-1$.
Now, 16 cases remain.
\item
Suppose that $\varepsilon_{12}=\varepsilon_{23}=1$,
all other $\varepsilon_{ij}=-1$.
Then the action of $\sigma_3$ on $z_i$ given in (\ref{3minus}) is
$$\sigma_3: z_1 \mapsto z_1, \, z_2 \mapsto -z_2, \,
z_3 \mapsto \frac{c}{z_3}(z_1^2-4az_2^2+16ab),$$
whose numerator is a linear polynomial of $z_2^2$ over $K(z_1)$,
so that the Noether problem is affirmative.

The same holds also for 6 cases
which are obtained by the permutation of suffices from the condition
$\varepsilon_{12}=\varepsilon_{23}=1$, other $\varepsilon_{ij}=-1$.
Now, 10 cases remain.
\item \label{e12e21plus}
Suppose that $\varepsilon_{12}=\varepsilon_{21}=1$, all other $\varepsilon_{ij}=-1$.
Then the actions of $\sigma_2$ and $\sigma_3$ on $y_i$ given in (\ref{2D2}) are
$$\begin{cases}
\sigma_2: y_1 \mapsto y_1, \, y_2 \mapsto \frac{b}{y_2}, \,
y_3 \mapsto -y_3 & \\
\sigma_3: y_1 \mapsto -y_1, \, y_2 \mapsto -y_2, \,
y_3 \mapsto -\frac{c(y_1^2-4a)}{y_3}. &
\end{cases}$$
So $K(x_1,x_2,x_3)^{<\sigma_1,\sigma_2>}=K(y_1,y_2,y_3)^{<\sigma_2>}
=K(\zeta_1,\zeta_2,\zeta_3)$, where
$\zeta_1=y_1$, $\zeta_2=y_2+\frac{b}{y_2}$,
$\zeta_3=\big(y_2-\frac{b}{y_2}\big)y_3$,
on which $\sigma_3$ acts as
$$\zeta_1 \mapsto -\zeta_1, \, \zeta_2 \mapsto -\zeta_2, \,
\zeta_3 \mapsto \left(y_2-\frac{b}{y_2}\right)
\frac{c(y_1^2-4a)}{y_3}
=\frac{c(\zeta_2^2-4b)(\zeta_1^2-4a)}{\zeta_3}.$$
So the Noether problem for $G$ is reduced to the Noether problem
for the following $\rho$.
$$\rho: x \mapsto -x, \, y \mapsto -y, \,
z \mapsto \frac{c(x^2-a)(y^2-b)}{z}$$
This problem depends only on the $K^{\times2}$-cosets of $a,b,c$.
We shall discuss this problem in the next sectioin.

Similar situation occurs for $\varepsilon_{13}=\varepsilon_{31}=1$
or $\varepsilon_{23}=\varepsilon_{32}=1$.
Now 7 cases remain.
\item \label{5minus}
Suppose that $\varepsilon_{32}=1$, all other $\varepsilon_{ij}=-1$.
Then the actions of $\sigma_2$ and $\sigma_3$
on $y_i$ given in (\ref{2D2}) are
$$\begin{cases}
\sigma_2: y_1 \mapsto -y_1, \,
y_2 \mapsto -\frac{b(y_1^2-4a)}{y_2}, \, y_3 \mapsto y_3 & \\
\sigma_3: y_1 \mapsto -y_1, \, y_2 \mapsto -y_2, \,
y_3 \mapsto -\frac{c(y_1^2-4a)}{y_3}. &
\end{cases}$$
So $K(x_1,x_2,x_3)^{<\sigma_1,\sigma_2>}=K(y_1,y_2,y_3)^{<\sigma_2>}
=K(u_1,u_2,u_3)$
where $u_1=\frac{1}{y_1}\big\{y_2+\frac{b(y_1^2-4a)}{y_2}\big\}$,
$u_2=y_2-\frac{b(y_1^2-4a)}{y_2}$, $u_3=y_3$,
on which $\sigma_3$ acts as
$$u_1 \mapsto u_1, \, u_2 \mapsto -u_2, \,
u_3 \mapsto -\frac{c(y_1^2-4a)}{u_3}=-\frac{c}{u_3}
\frac{u_2^2-4au_1^2}{u_1^2-4b},$$
whose numerator is a linear polynomial of $u_2^2$ over $K(u_1)$,
so that the Noether problem is affirmative.

The same holds for 6 cases in which just one of $\varepsilon_{ij}$ is $-1$.
Now, only one case remains.
\item \label{allminus}
Suppose that all $\varepsilon_{ij}$ are $-1$.
Then the action of $\sigma_3$ on $u_i$ given in (\ref{5minus}) is
$$\sigma_3: u_1 \mapsto -u_1, \, u_2 \mapsto u_2, \,
u_3 \mapsto -\frac{c}{u_3}\frac{u_2^2-4au_1^2}{u_1^2-4b}.$$
Thus in this case, the Noether problem for $G$
is reduced to the Noether problem for the following $\rho^\prime$.
$$\rho^\prime: x \mapsto x, \, y \mapsto -y, \,
z \mapsto -\frac{c}{z}(x^2-ay^2)(y^2-b)$$
This problem depends only on the $K^{\times2}$-cosets of $a,b,c$.
\end{sub}

The cases i) and ii) in Theorem \ref{3331} are reduced to
the Noether problems for $\rho$ and $\rho^\prime$
in (\ref{e12e21plus}) and (\ref{allminus}).
We shall consider them in the next section.
\end{pf}

\section{Problems $R_1(a,b,c)$ and $R_2(a,b,c)$.}

Consider the following two problems.

$R_1(a,b,c)$: The Noether problem for
$$\rho: x \mapsto -x, \, y \mapsto -y, \,
z \mapsto \frac{c(x^2-a)(y^2-b)}{z}$$.

$R_2(a,b,c)$: The Noether problem for
$$\rho^\prime: x \mapsto x, \, y \mapsto -y, \,
z \mapsto -\frac{c}{z}(x^2-ay^2)(y^2-b)$$.

Then, Theorem \ref{3331} is a result of the following Theorem.

\begin{thm}
$R_1(a,b,c)$ (resp. $R_2(a,b,c)$) is affirmative if and only if
$[K(\sqrt{a},\sqrt{b},\sqrt{c}):K] \leq 4$.
\end{thm}

\begin{pf}
First note that $R_1(a,b,c)$ is equivalent to $R_2(a,b,ac)$,
because $\rho$ is rewritten using $\xi=\frac{y}{x}$ and $\zeta=\xi z$ as
$$\rho: \xi\mapsto \xi, \, y \mapsto -y, \,
\zeta \mapsto \frac{c}{\zeta}(y^2-a\xi^2)(y^2-b)$$

The problem $R_1(a,b,c)$ does not change by interchanging $a$ and $b$,
because of the symmetricity of $\rho$.
The problem $R_2(a,b,c)$ does not change by any permutation of $a,b,c$,
because of  the symmetricity of the Noether problem
of the group (3,3,3,1) with all $\varepsilon_{ij}=-1$.

If $a=1$, then $R_1(a,b,c)$ is affirmative,
because putting $z^\prime=\frac{z}{x+1}$,
we have $\rho: z^\prime \mapsto -\frac{c}{z^\prime}(y^2-b)$.

From the symmetricity of $R_1$, it is affirmative also when $b=1$,
so that $R_2(a,b,c)$ is affirmative when $a=1$ or $b=1$.
From the symmetricity of $R_2$, it is affirmative also when $c=1$,
which in turn implies that $R_1$ is affirmative when $ac=1$.
It is affirmative also when $bc=1$,
which in  turn implies that $R_2$ is affirmatve when $abc=1$.
Thus we get the following result.

\begin{sub}
\item[(A)]
$R_1$ is affirmative when one of $a$, $b$, $ac$, $bc$ belongs to $K^{\times2}$.
$R_2$ is affirmative when one of $a$, $b$, $c$, $abc$ belongs to $K^{\times2}$.
\end{sub}

Next, consider the Noether problem of the group (3,3,3,1)
with $\varepsilon_{12}=\varepsilon_{21}=1$, all other $\varepsilon_{ij}=-1$.
If $c=1$, putting $\xi_3=\frac{x_3-1}{x_3+1}$,
we have $\sigma_3: \xi_3 \mapsto -\xi_3$,
so that $K(x_1,x_2,x_3)^{<\sigma_3>}=K(\eta_1,\eta_2,\eta_3)$
where $\eta_1=\frac{x_1}{\xi_3}$, $\eta_2=x_2\xi_3$, $\eta_3=x_1 \xi_3$
on which $\sigma_1$ and $\sigma_2$ act as
$$\begin{cases}
\sigma_1: \eta_1 \mapsto \frac{a}{\eta_1}, \,
\eta_2 \mapsto \frac{\eta_1\eta_2}{\eta_3}, \,
\eta_3 \mapsto \frac{a}{\eta_3} & \\
\sigma_2: \eta_1 \leftrightarrow \eta_3, \,
\eta_2 \mapsto \frac{b}{\eta_2} &
\end{cases}$$
Therefore $K(x_1,x_2,x_3)^{<\sigma_1,\sigma_3>}
=K(\eta_1,\eta_2,\eta_3)^{<\sigma_1>}=K(v_1,v_2,v_3)$
where $v_1=(\eta_1+\eta_3)\big/(\eta_1\eta_3+a)$,
$v_2=(\eta_1-\eta_3)\big/(\eta_1\eta_3-a)$,
$v_3=(\eta_1/\eta_3+1)\eta_2$
on which $\sigma_2$ acts as
$$\sigma_2: v_1 \mapsto v_1, \, v_2 \mapsto -v_2, \,
v_3 \mapsto \left(\frac{\eta_3}{\eta_1}+1\right)\frac{b}{\eta_2}
=\frac{b}{v_3}\frac{4v_1^2(av_2^2-1)}{v_2^2-v_1^2}$$
so that
$v_3^\prime:=(v_1+v_2)v_3 \mapsto
-4bv_1^2(av_2^2-1)/v_3^\prime$.

This means that $K(x_1,x_2,x_3)^G$ is rational,
so that $R_1(a,b,c)$ is affirmative when $c=1$.
This implies that $R_2(a,b,c)$ is affirmative when $ac=1$,
so that from the symmetricity it is affirmative also when $bc=1$ or $ab=1$,
which in turn means that $R_1$ is affirmative when $abc=1$ or $ab=1$.
Thus we get the following result.

\begin{sub}
\item[(B)]
$R_1$ is affirmative when one of $c$, $ab$, $abc$ belongs to $K^{\times2}$.
$R_2$ is affirmative when one of $ac$, $bc$, $ab$ belongs to $K^{\times2}$.
\end{sub}

Combining the results (A) and (B),
we see that $R_1$ (resp. $R_2$) is affirmative
when $[K(\sqrt{a},\sqrt{b},\sqrt{c}):K] \leq 4$.

Finally, consider the case $[L:K]=8$,
where $L=K(\alpha,\beta,\gamma)$ with
$\alpha=\sqrt{a}$, $\beta=\sqrt{b}$, $\gamma=\sqrt{c}$.
We shall prove the negativity of the Noether problem by the parity test.
For the action of $\rho$ defining $R_1$,
put $\zeta=\frac{z-\gamma(x+\alpha)(y+\beta)}{z+\gamma(x+\alpha)(y+\beta)}$,
then $\rho$ maps $\zeta$ to $-\zeta$.
Thus $L(x,y,z)^{<\rho>}=L(u,v,w)$
where $u=x^2$, $v=xy$, $w=x\zeta$.

$\mathfrak{G}=Gal(L/K)$ is isomorphic to $C_2 \times C_2 \times C_2$,
with generators $\tau_1: \alpha \mapsto -\alpha$, $\beta,\gamma$ invariant,
$\tau_2: \beta \mapsto -\beta$, $\alpha,\gamma$ invariant and
$\tau_3: \gamma \mapsto -\gamma$, $\alpha,\beta$ invariant.

Evidently $u,v$ are $\mathfrak{G}$-invariant.
The action of $\mathfrak{G}$ on $w$ is as follows.
$$\begin{cases}
{\displaystyle
\tau_3: \zeta \mapsto \frac{1}{\zeta},
\mbox{ so that } w \mapsto \frac{u}{w}} & \\
{\displaystyle \tau_1: \zeta \mapsto \frac{x\zeta+\alpha}{\alpha\zeta+x},
\mbox{ so that } w \mapsto \frac{u(w+\alpha)}{\alpha w+u}} & \\
{\displaystyle \tau_2: \zeta \mapsto \frac{y\zeta+\beta}{\beta\zeta+y},
\mbox{ so that } w \mapsto \frac{vw+\beta u}{\beta w+v}} &
\end{cases}$$
The smallest $\mathfrak{G}$-module containing $u$, $v$, $w$ is the following $M$ of rank 11.
The generators are $u$, $v$, $u-a$, $v^2-bu$, $w$, $\alpha w+u$, $\beta w+v$,
$w+\alpha$, $vw+\beta u$, $\beta u(w+\alpha)+v(\alpha w+u)$
and $vw+\beta u+\alpha(\beta w+v)$.
The action of $\mathfrak{G}$ is represented by matrices as follows.
$$m(\tau)=\begin{pmatrix}
I_4 & {\bf 0} \\
A(\tau) & B(\tau)
\end{pmatrix}$$
$$\begin{matrix}
A(\tau_1)=\begin{pmatrix}
 1 & 0 & 0 & 0 \\
 1 & 0 & 1 & 0 \\
 0 & 0 & 0 & 0 \\
 0 & 0 & 1 & 0 \\
 1 & 0 & 0 & 0 \\
 1 & 0 & 1 & 0 \\
 0 & 0 & 1 & 0
\end{pmatrix}, \, &
B(\tau_1)=\begin{pmatrix}
 0 &-1 & 0 & 1 & 0 & 0 & 0 \\
 0 &-1 & 0 & 0 & 0 & 0 & 0 \\
 0 &-1 & 0 & 0 & 0 & 1 & 0 \\
 1 &-1 & 0 & 0 & 0 & 0 & 0 \\
 0 &-1 & 0 & 0 & 0 & 0 & 1 \\
 0 &-1 & 1 & 0 & 0 & 0 & 0 \\
 0 &-1 & 0 & 0 & 1 & 0 & 0
\end{pmatrix} \\
A(\tau_2)=\begin{pmatrix}
 0 & 0 & 0 & 0 \\
 0 & 0 & 0 & 0 \\
 0 & 0 & 0 & 1 \\
 0 & 0 & 0 & 0 \\
 0 & 0 & 0 & 1 \\
 0 & 0 & 0 & 1 \\
 0 & 0 & 0 & 1
\end{pmatrix}, \, &
B(\tau_2)=\begin{pmatrix}
 0 & 0 &-1 & 0 & 1 & 0 & 0 \\
 0 & 0 &-1 & 0 & 0 & 1 & 0 \\
 0 & 0 &-1 & 0 & 0 & 0 & 0 \\
 0 & 0 &-1 & 0 & 0 & 0 & 1 \\
 1 & 0 &-1 & 0 & 0 & 0 & 0 \\
 0 & 1 &-1 & 0 & 0 & 0 & 0 \\
 0 & 0 &-1 & 1 & 0 & 0 & 0
\end{pmatrix} \\
A(\tau_3)=\begin{pmatrix}
 1 & 0 & 0 & 0 \\
 1 & 0 & 0 & 0 \\
 0 & 0 & 0 & 0 \\
 0 & 0 & 0 & 0 \\
 1 & 0 & 0 & 0 \\
 1 & 0 & 0 & 0 \\
 0 & 0 & 0 & 0
\end{pmatrix}, \, &
B(\tau_3)=\begin{pmatrix}
-1 & 0 & 0 & 0 & 0 & 0 & 0 \\
-1 & 0 & 0 & 1 & 0 & 0 & 0 \\
-1 & 0 & 0 & 0 & 1 & 0 & 0 \\
-1 & 1 & 0 & 0 & 0 & 0 & 0 \\
-1 & 0 & 1 & 0 & 0 & 0 & 0 \\
-1 & 0 & 0 & 0 & 0 & 0 & 1 \\
-1 & 0 & 0 & 0 & 0 & 1 & 0
\end{pmatrix}
\end{matrix}$$
From this, we can calculate the cohomology group $H^1$ and $\hat{H}^{-1}$.
As a result, we get $\hat{H}^{-1}(\mathfrak{G},M)=0$
so that the non-vanishing cohomology test does not work.

On the other hand, we have $H^1(\mathfrak{H},M) \not= 0$
for $|\mathfrak{H}| \geq 2$.
However, if we extend $M$ to a module $M^\prime$ of rank 12
by adding $w^2-u$ as the twelfth generator,
then we get $H^1(\mathfrak{H},M^\prime)=0$
for any subgroup $\mathfrak{H}$ of $\mathfrak{G}$.
So we shall apply the parity test for this $M^\prime$.

The action of $\mathfrak{G}$ on $M^\prime$ is represented as matrices
as follows.
$$m^\prime(\tau)=\begin{pmatrix}
m(\tau) & {\bf 0} \\
C(\tau) & 1
\end{pmatrix}$$
$$\begin{array}{ccccccccccccccc}
C(\tau_1) &=& ( &
 1 & 0 & 1 & 0 & 0 &-2 & 0 & 0 & 0 & 0 & 0
& ) \\
C(\tau_2) &=& ( &
 0 & 0 & 0 & 1 & 0 & 0 &-2 & 0 & 0 & 0 & 0
& ) \\
C(\tau_3) &=& ( &
 1 & 0 & 0 & 0 &-2 & 0 & 0 & 0 & 0 & 0 & 0
& )
\end{array}$$
If $K(x,y,z)^{<\rho>}$ is rational,
then $M^\prime$ is a direct product factor of a permutation module $Q$.
Let $i$ be the injection $M^\prime \rightarrow Q$
and $p$ be the projection $Q \rightarrow M^\prime$,
then we have $p \circ i=id_{M^\prime}$.

Since $u$ is $\mathfrak{G}$-invariant,
$i(u)$ has a constant exponent on every transitive part
of the permutations induced by $\mathfrak{G}$.
$Q$ is a direct sum of transitive parts $X$,
and $i(u)=\sum_X q_{{}_X}$ where each $q_{{}_X}$ is in the form of
$m \sum_{\overline{\tau} \in \mathfrak{G}/\mathfrak{G}^\prime} \beta^{\overline{\tau}}$.
Here $m$ is the constant exponent of $i(u)$ on X,
$\mathfrak{G}^\prime$ is the stabilizer of $X$
and $\beta$ is some fixed element of $X$.
Since $p \circ i(u)=\sum_X p(q_{{}_X})$,
a contradiction occurs if the exponent of $u$ of each $p(q_{{}_X})$ is even.

Since $\tau_3$ maps $w$ to $\frac{u}{w}$,
if $\tau_3 \in \mathfrak{G}^\prime$,
then we have $n=m-n$ where $n$ is the exponent of $w$ on $X$,
so $m$ must be even.
Thus we get $p(q_{{}_X}) \in 2M^\prime$.

Suppose that $\tau_3 \not\in \mathfrak{G}^\prime$.
Let $\mathfrak{H}$ be a subgroup of order 4
which contains $\mathfrak{G}^\prime$
but does not contain $\tau_3$.
It is sufficient to assume $m=1$,
then $q_{{}_X}=\varphi^\tau+\varphi$
where $\varphi=\sum_{\overline{\tau} \in \mathfrak{H}/\mathfrak{G}^\prime} \beta^{\overline{\tau}}$,
and $\tau \not\in \mathfrak{H}$.

If $\mathfrak{H} \not= <\tau_2, \tau_1\tau_3>$,
then one of $\tau_2$, $\tau_1\tau_3$, $\tau_1\tau_2\tau_3$
does not belong to $\mathfrak{H}$,
so that $p(q_{{}_X})=\big(m^\prime(\tau)+1\big)p(\varphi) \in \im(m^\prime(\tau)+1\big)$,
but the first column of $m^\prime(\tau)+1$ is
$\begin{pmatrix} 2 \\ 0 \\ \vdots \\ 0 \end{pmatrix}$
when $\tau$ is one of $\tau_2$, $\tau_1\tau_3$, $\tau_1\tau_2\tau_3$,
so that $p(q_{{}_X})$ has an even exponent on $u$.

If $\mathfrak{H}=<\tau_2,\tau_1\tau_3>$,
then $p(q_{{}_X})=\big(m^\prime(\tau_3)+1\big)p(\varphi)
\in \big(m^\prime(\tau_3)+1\big)(M^{\prime\mathfrak{H}})$
since $\varphi$ is $\mathfrak{H}$-invariant.

Denoting the group operation of $M^\prime$ additively,
let $e_1 \sim e_{12}$ be the basis of $M^\prime$ defined before.
Then $M^{\prime\mathfrak{H}}$ is a $\mathbb{Z}$-module of rank 6
and generated by
$f_1=e_5+e_6+e_9+e_{10}-2e_{12}$ and $f_2=e_7+e_8+e_{11}-2e_{12}$
besides $e_1$, $e_2$, $e_3$, $e_4$.

Regarding $e_1 \sim e_4$ and $f_1$, $f_2$
as the basis of $M^{\prime\mathfrak{H}}$,
$\tau_3+1$ is represented as a (6,12) matrix.
Then the first column becomes
$\begin{pmatrix} 2 \\ 0 \\ 0 \\ 0 \\ 2 \\ -2 \end{pmatrix}$,
so that $p(q_{{}_X})$ has an even exponent on $u$.
Anyway $p \circ i(u)$ has an even exponent on $u$,
which is a contradiction to $p \circ i=id_{M^\prime}$.
Thus $R_1(a,b,c)$ is negative when
$[K(\sqrt{a},\sqrt{b},\sqrt{c}):K]=8$.
\end{pf}

\end{document}